\documentclass{amsart}

\newcommand{\divides}{|}

\newcommand{\fix}{\operatorname{\mathsf F}}

\newcommand{\orbit}{\operatorname{\mathsf O}}
\newcommand{\mertens}{\operatorname{\mathsf M}}

\usepackage{amsmath,amsfonts,amsthm,amssymb,psfrag,graphicx}

\DeclareMathOperator\ord{ord}

\newtheorem{theorem}{Theorem}
\newtheorem{lemma}[theorem]{Lemma}

\theoremstyle{definition}
\newtheorem{example}{Example}

\renewcommand{\le}{\leqslant}
\renewcommand{\ge}{\geqslant}
\renewcommand{\epsilon}{\varepsilon}

\newcommand{\sg}{\mathcal L}

\def\bigo{\operatorname{O}}    
\def\littleo{\operatorname{o}} 

\def\lowerboundpi{1}
\def\upperboundpi{2}
\def\mertensnilpotent{3}
\def\mertensabelian{4}
\def\mobiusboundnilpotent{5}
\def\firstabelianupperbound{6}
\def\secondabelianupperbound{7}
\def\dusautoygrunewaldbound{8}
\def\abeliansnbound{9}

\title{Orbit-counting for nilpotent
group shifts}
\author{Richard Miles}
\author{Thomas Ward}
\date{\today}

\address{School of Mathematics, University of East Anglia,
Norwich, NR4 7TJ, UK}
\thanks{
We thank Johannes Siemons and Shaun Stevens for their
suggestions.
This research was supported by E.P.S.R.C.
grant EP/C015754/1.}

\begin{document}

\begin{abstract}
We study the asymptotic behaviour of the orbit-counting
function and a dynamical Mertens' theorem for
the full~$G$-shift for a finitely-generated torsion-free
nilpotent group~$G$. Using bounds for the M{\"o}bius function on
the lattice of subgroups of finite index and
known subgroup growth estimates, we find a
single asymptotic of the shape
\[
\sum_{\vert\tau\vert\le N}\frac{1}{e^{h\vert\tau\vert}}\sim
CN^{\alpha}(\log N)^{\beta}
\]
where~$\vert\tau\vert$ is the cardinality of the
finite orbit~$\tau$. For the usual orbit-counting
function we find upper and lower bounds together
with numerical evidence to suggest that
for actions of non-cyclic groups there is no
single asymptotic in terms of elementary functions.
\end{abstract}
\maketitle

\section{Introduction}\label{sec:introduction}

There is a well-understood analogy between the
classical prime number theorem and Mertens' theorem on
the one hand, and closed orbits of dynamical systems
on the other.
For reasons that will become clear,
we will think of an invertible map~$T:X\to X$
as defining a~$\mathbb Z$-action, also written~$T$,
by bijections of~$X$. For any subgroup~$L\leqslant\mathbb Z$,
write
\begin{equation}\label{def:fixL}
\fix_T(L)=\vert\{x\in X\mid
T^{\ell}(x)=x\mbox{ for all }\ell\in L\}\vert
\end{equation}
for the number of~$L$-periodic
points under the action,
and~$[L]$ for the index of~$L$ in~$\mathbb Z$. Write
\begin{equation}\label{def:orbitL}
\orbit_T(L)=
\frac{1}{[L]}
\vert
\{
x\in X\mid
T^{\ell}(x)=x\Leftrightarrow\ell\in L\}
\vert
\end{equation}
for the number of~$L$-orbits.
Noting that~$L$ is determined by
the elements of an~$L$-orbit,
write~$\vert\tau\vert=[L]$
if~$\tau$ is an~$L$-orbit, and call a set~$\tau$
a
\emph{closed orbit} if it is an~$L$-orbit for some~$L$
of finite index.

Since~$\fix_T(L)=\displaystyle\sum_{L'\geqslant L}
[L']\orbit_T(L')$, M{\"o}bius inversion
shows that
\begin{equation}\label{def:mobiusinversion1}
\orbit_T(L)=
\sum_{L'\geqslant L}
\mu\left(
\left\vert L'/L\right\vert
\right)\fix_T(L').
\end{equation}
If~$T$ is hyperbolic with respect to some metric
structure on~$X$, and~$T$ has topological entropy~$h$,
then the dynamical analogue of the prime
number theorem takes the form
\[
\pi_T(N)\sim
\frac{e^{(N+1)h}}{N}
\]
where
\begin{equation}\label{whydidithavetostop}
\pi_T(N)=\left\vert
\{\tau\mid
\tau\mbox{ is a closed orbit with }
\vert\tau\vert\leqslant N\}
\right\vert
\end{equation}
and the dynamical analogue of Mertens' theorem is
\[
\mertens_T(N)\sim\log N+C+\littleo(1)
\]
where
\begin{equation}\label{youveblownitallskyhigh}
\mertens_T(N)=\sum_{\vert\tau\vert\leqslant N}\frac{1}{e^{h\vert\tau\vert}}.
\end{equation}
These are simplified statements of
results due to Parry and
Pollicott~\cite{MR727704}, Sharp~\cite{MR1139566}
and others, motivated by results
of Margulis~\cite{MR2035655} on orbits
of geodesic flows.

Our purpose here is to explore possible extensions of
these results to the setting of other group actions.
For algebraic~$\mathbb Z^d$-actions, the growth
rate of periodic points has been studied by
Lind, Schmidt, Ward~\cite[Sec.~7]{MR1062797}, and some properties
of a natural dynamical zeta function have been
studied by Lind~\cite{MR1411232}.
Despite this, the behaviour of closed orbits in these
systems presents many difficulties: in particular,
the growth rate of periodic points is only really
understood for periods with good geometrical properties
(corresponding to subgroups~$L$ that are
generated by vectors of length comparable to~$[L]^{1/d}$),
and the dynamical zeta function is already a highly
non-trivial object for the trivial action on a point.
In order to find analogues of the prime number
theorem and Mertens' theorem we need to take account
of all orbits, no matter how awkwardly shaped.
Because of the difficulties presented by this,
and in order to expose some distinctive features of
genuinely non-cyclic actions, we restrict
attention to full shifts, which are defined
as follows. For any group~$G$ and finite alphabet~$B$,
define the full~$G$-shift~$T:B^G\to B^G$
by
\[
(
T_g(x))_h=x_{gh}.
\]
It is clear that if~$L\leqslant G$ is a subgroup
of finite index, then
\[
\fix_T(L)=b^{[L]}
\]
where~$b=\vert B\vert$, and
the topological entropy of the full shift is equal to~$\log b$.

While our results are stated for this special kind
of system, one of our aims is to highlight the
way in which subgroup growth
is inherently linked to orbit-growth properties.
This is
also
true for single transformations but is much less
noticeable because~$\mathbb Z$ happens to have
exactly one subgroup of each index.

\section{Closed orbits of group actions}

Now consider an action~$T$ of a
finitely-generated group~$G$ by
bijections of a set~$X$.
Let~$\sg=\sg(G)$ denote the set
of subgroups of~$G$ with finite index, and write~$[L]$
for the index of~$L$ in~$G$.
Define~$\fix_T(L)$
and~$\orbit_T(L)$ using~\eqref{def:fixL} and~\eqref{def:orbitL}
for any~$L\in\sg$. In order to write down the
analogue of~\eqref{def:mobiusinversion1},
some
terminology from combinatorics is needed.
The partial order~$\preccurlyeq$ defined by~$L'\preccurlyeq L$ if and only
if~$L\le L'$ makes~$\sg$ into a locally finite
poset, which therefore has a M{\"o}bius function
defined on intervals as follows
(see Stanley~\cite[Sec.~3.7]{MR1442260} for
the details and examples; for convenience we use
the relations~$\le,\ge$ rather than~$\succcurlyeq,\preccurlyeq$
throughout).
Define a function~$\mu$ on pairs~$(L',L)\in\sg^2$
with~$L\leqslant L'$
by the properties
\[
\mu(L,L)=1\mbox{ for all }L\in\sg
\]
and
\[
\mu(L',L)=-\sum_{L< L''\leqslant L'}
\mu(L',L'')\mbox{ for all }L< L';
\]
the symbol~$\mu$ will always denote the M{\"o}bius function
corresponding to~$G$, for a different group~$H$ we write~$\mu_H$.
For the case~$G=\mathbb Z$ dealt with in Section~\ref{sec:introduction},
the set~$\sg$ is in one-to-one correspondence
with~$\mathbb N\cup\{0\}$ via the
correspondence~$d\to d\mathbb Z$. Under
this correspondence~$d\mathbb Z\leqslant d'\mathbb Z$ if
and only if~$d'\divides d$, and~$\mu(d'\mathbb Z,
d\mathbb Z)$ is just~$\mu(d/d')$, the classical M{\"o}bius
function. In order to understand one of the
difficulties that arises in passing from~$\mathbb Z$-actions
to
more general group actions, notice that
one of the key estimates behind the results
in~\cite{emsw}
is~$\vert\mu\vert\leqslant 1$
on~$\mathbb Z$, while on larger groups there is no
reason to expect~$\mu$ to be bounded.

For any group~$G$, define the rank of~$G$, written~$r(G)$,
to be the supremum over all
finitely-generated subgroups~$H$ of
the minimal cardinality of a generating set
for~$H$. 
The analogues of the prime number theorem
and Mertens' theorem in our
setting take the following form.

\begin{theorem}\label{lucyintheskywithdiamonds}
Let~$G$ be a finitely generated torsion-free
nilpotent group, and let~$T$ be the full~$G$-shift
on an alphabet with~$b$ symbols.
\begin{enumerate}
\item
For any~$\epsilon>0$ there are positive constants~$C_{\lowerboundpi}$
and~$C_{\upperboundpi}$ for
which
\begin{equation*}
C_{\lowerboundpi}N^{d(G)-2}\le\frac{\pi_T(N)}{b^N}\le C_{\upperboundpi}N^{r(G)-1}
\end{equation*}
where~$r(G)$ is the rank of~$G$ and~$d(G)$
is the rank of the abelianization of~$G$.
In the case~$G=\mathbb Z^d,d\ge2$,
\begin{equation*}
C_{\lowerboundpi}N^{d-2}\le\frac{\pi_T(N)}{b^N}\le C_{\upperboundpi}N^{d-2}(\log N)^{d-1}.
\end{equation*}
\item 
There is a positive constant~$C_{\mertensnilpotent}$, and there are
non-negative constants~$\alpha\in\mathbb Q$ and~$\beta\in\mathbb Z$,
with
\[
\mertens_T(N)\sim C_{\mertensnilpotent}N^{\alpha}(\log N)^{\beta},
\]
and in the case~$G=\mathbb Z^d,d\ge2$, there
is a positive constant~$C_{\mertensabelian}$ with
\[
\mertens_T(N)\sim C_{\mertensabelian}N^{d-1}.
\]
\end{enumerate}
\end{theorem}

\begin{proof}
We begin by assembling some bounds for~$a_n(G)$,
the number of subgroups of~$G$
with index~$n$. By
Grunewald, Segal and Smith~\cite[Prop.~1.1]{MR943928} (this is
also shown in a different setting by
Lind~\cite[Prop.~4.2]{MR1411232}) we have
\begin{equation}\label{itsthreeoclockinthemorning}
a_n(\mathbb Z^d)=\sum_{k\divides n}a_{n/k}(\mathbb Z^{d-1})k^{d-1},
\end{equation}
so, in particular,
\begin{equation}\label{boundsonanforzd}
a_n(\mathbb Z^d)\ge n^{d-1}\mbox{ for all }n,d\ge1.
\end{equation}
It follows that
\begin{equation}\label{lowerboundforannilpotent}
a_n(G)\ge n^{d(G)-1}.
\end{equation}
From~\cite[Lemma~1.4.1]{MR1978431},
\begin{equation}\label{upperboundforannilpotent}
a_n(G)<n^{r(G)}\mbox{ for }n>1,
\end{equation}
so~$s_n(G)=\sum_{k=1}^{n}a_k(G)<n^{r(g)+1}$.

Notice that
\begin{eqnarray}
\pi_T(N)=\sum_{\vert\tau\vert\le N}1&=&\sum_{[L]\le N}
\orbit_T(L)\nonumber\\
&=&
\sum_{[L]\le N}\frac{1}{[L]}\sum_{L'\ge L}
\fix_T(L')\mu(L',L)\nonumber\\
&=&
\sum_{[L]\le N}\frac{1}{[L]}\fix_T(L)+
\underbrace{
\sum_{[L]\le N}\frac{1}{[L]}
\sum_{L'>L}\fix_T(L')\mu(L',L)}_{\Sigma_N}.\label{bytellingmeyourlie}
\end{eqnarray}
We wish to show that~$\Sigma_N/b^N\to0$ as~$N\to\infty$.
Clearly
\[
F_T(L')\le b^{[L]/2}\mbox{ for all }L<L'.
\]

\begin{lemma}\label{takeiteasy}
There exists a constant~$C_{\mobiusboundnilpotent}$ 
such that~$\mu(L',L)\le e^{C_{\mobiusboundnilpotent}(\log[L])^2}$
for~$L\le L'$.
\end{lemma}

\begin{proof}
If~$L$ is not normal in~$L'$, then (since
passing to the normal closure is a closure
operator on~$\sg$), Crapo's Theorem (see~\cite[Th.~1]{MR0245483})
shows that~$\mu(L',L)=0$.

If~$L$ is normal in~$L'$, then~$\mu(L',L)=
\mu_{L'/L}(L'/L,\{0\})$, and by Kratzer and
Th{\'e}venaz~\cite[Prop.~2.4]{MR761806},~$\mu_{L'/L}(L'/L,\{0\})=
0$ if~$L'/L$ is not a product of elementary
abelian groups, and is bounded by
\[
\prod_{i=1}^{r}
p_i^{n_i(n_i-1)/2}
\]
if~$L'/L\cong\prod_{i=1}^{r}(\mathbb Z/p_i\mathbb Z)^{n_i}.$
In this case, writing~$\vert L'/L\vert=k$ and~$v_i(\cdot)=
\ord_{p_i}(\cdot)$, it follows that
\begin{eqnarray*}
\left(
\mu_{L'/L}(L'/L,\{0\})\right)^2\le
\prod_{i=1}^{r}p_i^{v_i(k)^2}
=
\prod_{i=1}^{r}\left(p_i^{v_i(k)}\right)^{v_i(k)}
\le
\left(\prod_{i=1}^{r}p_i^{v_i(k)}\right)^{\log_2 k}
=k^{\log_2(k)},
\end{eqnarray*}
so
there is
a constant~$C_{\mobiusboundnilpotent}$ with
\begin{equation*}
\vert\mu_{L'/L}(L'/L,\{0\})\vert\le
e^{C_{\mobiusboundnilpotent}(\log[L])^2}.
\end{equation*}
\end{proof}
Therefore,~\eqref{bytellingmeyourlie} gives
\begin{eqnarray*}
\frac{\Sigma_N}{b^N}&\le&\frac{1}{b^N}
\sum_{[L]\le N}\frac{1}{[L]}
b^{[L]/2}e^{C_{\mobiusboundnilpotent}(\log[L])^2}s_{[L]}(G)\\
&\le&
\frac{1}{b^N}\sum_{n=1}^{N}
\frac{1}{n}b^{n/2}e^{C_{\mobiusboundnilpotent}(\log n)^2}a_n(G)s_n(G)\\
&\le&
\frac{1}{b^N}\sum_{n=1}^{N}
\frac{1}{n}
b^{n/2}e^{C_{\mobiusboundnilpotent}(\log n)^2}n^{2r(G)+1}
\quad\mbox{ by }\eqref{upperboundforannilpotent}\\
&\le&
{b^{-N/2}}N^{2r(G)+2}e^{C_{\mobiusboundnilpotent}(\log N)^2}
\rightarrow 0\mbox{ as }N\to\infty.
\end{eqnarray*}
It follows that
\begin{equation}\label{itsagirlmylord}
\frac{\pi_T(N)}{b^N}+\littleo(1)=
\frac{1}{b^N}\sum_{[L]\le N}\frac{1}{[L]}F_T(L)
=\frac{1}{b^N}\sum_{n=1}^{N}
\frac{a_n(G)}{n}b^n.
\end{equation}
By partial summation,
\begin{equation}\label{inaflatbedford}
\sum_{n=1}^{N}n^eb^n=\frac{b}{b-1}N^eb^N+\bigo\left(N^{e-1}b^N\right)
\end{equation}
for any~$e\ge1$ since~$b\ge2$.
Thus by~\eqref{lowerboundforannilpotent}
and~\eqref{upperboundforannilpotent}, part~(1) of the
theorem in the nilpotent case
follows from~\eqref{itsagirlmylord}.

Now consider the case~$G=\mathbb Z^d$.
The estimate~\eqref{boundsonanforzd}
in place of~\eqref{lowerboundforannilpotent} gives
the stated lower bound
by the same argument.
For the upper bound we
use the following lemma.

\begin{lemma}\label{igaveyoulove}
For~$n\ge2$,
\[
a_n(\mathbb Z^d)\le3^dn^{d-1}(\log n)^{d-1}.
\]
\end{lemma}

\begin{proof}
If~$d=2$, then~\eqref{itsthreeoclockinthemorning}
gives~$a_n(\mathbb Z^2)=\sigma(n)=\sum_{d\divides n}d$, and
a simple argument shows that~$\sigma(n)\le
3n\log(n)$ for~$n\ge2$.
Now assume the statement of the lemma, and
notice that by~\eqref{itsthreeoclockinthemorning}
we have
\begin{eqnarray*}
a_n(\mathbb Z^{d+1})&=&\sum_{k\divides n}a_{n/k}
(\mathbb Z^d)k^d\\
&\le&3^d\sum_{k\divides n}(n/k)^{d-1}
\left(\log(n/k)\right)^{d-1}k^d\mbox{ by
hypothesis}\\
&\le&3^dn^{d-1}(\log n)^{d-1}\sum_{k\divides n}
k\\
&\le&3^{d+1}n^d(\log n)^d\mbox{ for }n\ge 2.
\end{eqnarray*}
\end{proof}

Now there are positive constants~$C_{\firstabelianupperbound}$ and~$C_{\secondabelianupperbound}$ with
\begin{eqnarray*}
\frac{1}{b^N}\sum_{n=1}^{N}
\frac{a_n(\mathbb Z^d)}{n}b^n&\le&
C_{\firstabelianupperbound}\frac{1}{b^N}
\sum_{n=1}^{N}n^{d-2}(\log n)^{d-1}b^n\mbox{
by Lemma~\ref{igaveyoulove}}\\
&\le&
C_{\firstabelianupperbound}
\frac{1}{b^{N}}
(\log N)^{d-1}\sum_{n=1}^{N}n^{d-2}b^n\\
&\le&
C_{\secondabelianupperbound}N^{d-2}(\log N)^{d-1}\mbox{ by~\eqref{inaflatbedford}},
\end{eqnarray*}
giving the upper bound by~\eqref{itsagirlmylord}.

Turning to the analogue of Mertens' theorem, notice that
\begin{eqnarray*}
\mertens_T(N)=\sum_{\vert\tau\vert\le N}
\frac{1}{b^{\vert\tau\vert}}&=&\sum_{n=1}^{N}
\frac{1}{b^n}\sum_{[L]=n}\orbit_T(L)\\
&=&
\sum_{n=1}^N\frac{1}{b^n}\sum_{[L]=n}\frac{1}{n}
\sum_{L'\ge L}\mu(L',L)\fix_T(L')\\
&=&\sum_{n=1}^{N}\frac{a_n(G)}{n}+\underbrace{
\sum_{n=1}^{N}\frac{1}{nb^n}\sum_{[L]=n}
\sum_{L'>L}\mu(L',L)b^{[L']}
}_{\Delta_N},
\end{eqnarray*}
and by~\eqref{upperboundforannilpotent} and Lemma~\ref{takeiteasy}
we have
\[
\vert\Delta_N\vert\le\sum_{n=1}^{N}\frac{1}{nb^n}
n^{r(G)+1}b^{n/2}e^{C_{\mobiusboundnilpotent}(\log(n/2))^2}=\bigo(1).
\]
It follows that
\[
\mertens_T(N)=\sum_{n=1}^{N}\frac{a_n(G)}{n}+\bigo(1).
\]
Now by a deep theorem of
du Sautoy and Grunewald~\cite[Th.~1.1]{MR1815702},
there is a constant~$C_{\dusautoygrunewaldbound}$ and
there are
non-negative constants~$\gamma\in\mathbb Q$,~$\delta\in\mathbb Z$ with
\[
s_n(G)\sim C_{\dusautoygrunewaldbound}n^{\gamma}(\log n)^{\delta}.
\]
If~$\gamma\ge1$, then it follows by partial
summation that~$\mertens_T(N)$ has
the asymptotic stated in part~(2) with~$\alpha=\gamma-1$
and~$\beta=\delta$.
If~$\gamma<1$ then partial summation shows that~$\mertens_T(N)$
is bounded, so we may take~$\alpha=\beta=0$.
In the case~$G=\mathbb Z^d$,
the well-known relation
\[
\sum_{n=1}^{\infty}\frac{a_n(\mathbb Z^d)}{n^z}
=
\zeta(z)\zeta(z-1)\cdots\zeta(s-d+1)
\]
gives, via a Tauberian theorem (see~\cite{MR1815702}), an asymptotic
of the form
\[
s_n(\mathbb Z^d)\sim C_{\abeliansnbound}n^d,
\]
which gives the conclusion.
\end{proof}

\section{Examples}

There are two clear directions in which
Theorem~\ref{lucyintheskywithdiamonds}
is limited. It applies to a very special
class of dynamical systems, and
Examples~\ref{westoodinthewindycity}
and~\ref{thegipsyboyandi} indicate
some of the obstacles to extending it to
a broader class of systems.
Moreover, the estimates used are
na\"{\i}ve and no doubt more sophisticated
techniques could improve the bounds.

\begin{example}
It is not surprising that~$\mertens_T$ is less
sensitive to volatility in~$\orbit_T(L)$ than is~$\pi_T$;
in this example we
explain something about the extent of this difference
in the simplest non-cyclic case~$G=\mathbb Z^2$
with~$b=2$.
By Theorem~\ref{lucyintheskywithdiamonds},~$\mertens_T(N)\sim C_{\mertensnilpotent}N$.
On the other hand,
by~\eqref{itsthreeoclockinthemorning} we
have~$a_n(\mathbb Z^2)=\sigma(n)$, so
from~\eqref{itsagirlmylord}
we have
\begin{equation}\label{ithoughtthatwehadmadeittothetop}
\frac{\pi_T(N)}{2^N}=\frac{1}{2^N}\sum_{n=1}^{N}
\frac{\sigma(n)}{n}2^n+\littleo(1).
\end{equation}
The graph in Figure~\ref{java}
illustrates the asymptotic
in~\eqref{ithoughtthatwehadmadeittothetop} and
the erratic behaviour of~$\pi_T(N)$ as a 
function of~$N$.

\begin{figure}[htbp]
      \begin{center}
\psfrag{genuine}{\hskip1cm$\phi(N)$}
\psfrag{asymptotic}{\hskip1.4cm$\psi(N)$}
\psfrag{N}{$N$}
\psfrag{ 0}{$0$}
\psfrag{ 10}{$10$}
\psfrag{ 20}{$20$}
\psfrag{ 30}{$30$}
\psfrag{ 40}{$40$}
\psfrag{ 50}{$50$}
\psfrag{ 60}{$60$}
\psfrag{ 70}{$70$}
\psfrag{ 80}{$80$}
\psfrag{ 90}{$90$}
\psfrag{ 100}{$100$}
\psfrag{ 1}{$1$}
\psfrag{ 1.5}{$1.5$}
\psfrag{ 2}{$2$}
\psfrag{ 2.5}{$2.5$}
\psfrag{ 3}{$3$}
\psfrag{ 3.5}{$3.5$}
\psfrag{ 4}{$4$}
\psfrag{ 4.5}{$4.5$}
\scalebox{0.8}{\includegraphics{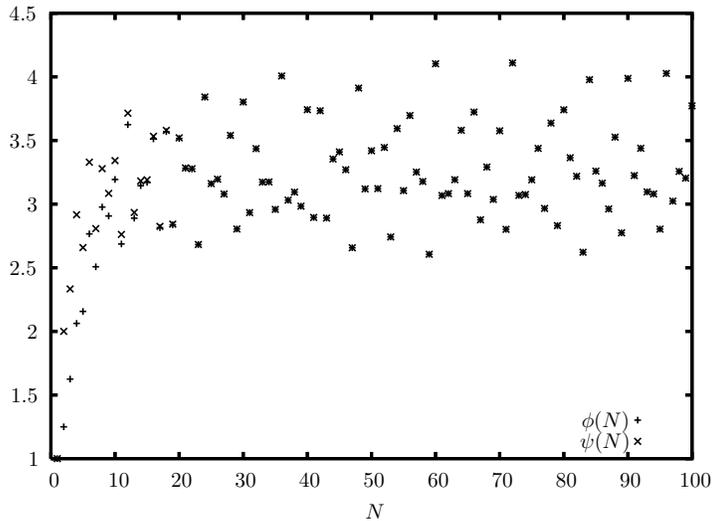}}
         \caption{\label{java}Plot of $\phi(N)=\displaystyle\frac{\pi_T(N)}{2^N}$ and
$\psi(N)=\displaystyle\frac{1}{2^N}\sum_{n=1}^N\frac{\sigma(n)}{n}2^n$.}
\end{center}
\end{figure}

\end{example}

\begin{example}
A simple example with~$\beta>0$ in
Theorem~\ref{lucyintheskywithdiamonds}(2)
is given by the discrete Heisenberg
group
\[
G=\left\{
\left(
\begin{matrix}1&a&b\\
0&1&c\\
0&0&1\end{matrix}
\right)\mid
a,b,c\in\mathbb Z
\right\}
\]
which has
\[
s_n(G)\sim\frac{\zeta(2)^2}{2\zeta(3)}
N^2\log N
\]
by Smith~\cite{smith}
(see also~\cite[Sect.~1]{MR1815702}).
It follows that
\begin{eqnarray*}
\mertens_T(N)&\sim&\frac{s_N(G)}{N}+
\int_1^N\log t{\thinspace\rm d} t\\
&\sim&\left(\frac{\zeta(2)^2}{2\zeta(3)}
+1\right)
N\log N.
\end{eqnarray*}
The bounds for~$\pi_T(N)/b^N$ are of course
much weaker,
\[
0<C_{\lowerboundpi}\le\frac{\pi_T(N)}{b^N}\le C_{\upperboundpi}N^{2}.
\]
\end{example}

A natural setting to seek results of this sort would
be expansive algebraic~$\mathbb Z^d$-actions, because
their growth rate of periodic points is
available from~\cite{MR1062797}.
However, this growth rate result concerns a
sequence of subgroups~$L_n$ with the property
that~$\mathsf{d}\left(0,L_n\setminus\{0\}
\right)\rightarrow\infty$
as~$n\to\infty$ (for the Euclidean
metric~$\mathsf{d}$ on~$\mathbb Z^d$).
As mentioned in the introduction,
genuine orbit-counting
results should include all orbits.
The examples above describe some of the
implications of erratic behavour in the
map~$n\mapsto a_n(G)$; the next two examples
concern two different ways in
which the map~$\fix:L\mapsto\fix_T(L)$ may be
erratic, in contrast to the full shift where~$\fix_T(L)$
is a regularly growing function of the
index~$[L]$ alone.

\begin{example}\label{westoodinthewindycity}
For expansive~$\mathbb Z^d$-actions without
any entropy assumption,~$\fix$ may be
sensitive to small changes in subgroups of
fixed shape.
For example,
in Ledrappier's~$\mathbb Z^2$-action (see~\cite{MR512106}) we
have
\[
\fix_T((2^k,0)\mathbb Z
\oplus(0,2^k)\mathbb Z)=1
\]
for all~$k\ge1$
and
\[
\fix_T((2^k-1,0)
\mathbb Z\oplus(0,2^k-1)\mathbb Z)\to\infty
\]
as~$k\to\infty$ (see~\cite[Ex.~3.3]{MR1165372}).
\end{example}

\begin{example}\label{thegipsyboyandi}
Expansive~$\mathbb Z^d$ actions with
completely positive entropy do have the
property that~$\fix_T(L)$ is not very 
sensitive to changes in~$[L]$ for a fixed
shape of subgroup. However,
they may still be sensitive to the
shape of~$L$ for a fixed index~$[L]$.
For example, the~$\mathbb Z^2$-action
on~$X^{\mathbb Z}$, where~$X$
is the Pontryagin dual of~${\mathbb Z[\frac{1}{2}]}$,
defined by~$\left(T_{(1,0)}(x)\right)_k=2x_k\pmod{1}$
and~$\left(T_{(0,1)}(x)\right)_k=x_{k+1}$
has
\[
\fix_T((n,0)\mathbb Z\oplus(0,1)\mathbb Z)=2^n-1,
\]
while
\[
\fix_T((1,0)\mathbb Z\oplus(0,n)\mathbb Z)=1
\]
for all~$n\ge1$.
\end{example}



\providecommand{\bysame}{\leavevmode\hbox to3em{\hrulefill}\thinspace}
\providecommand{\MR}{\relax\ifhmode\unskip\space\fi MR }
\providecommand{\MRhref}[2]{%
  \href{http://www.ams.org/mathscinet-getitem?mr=#1}{#2}
}
\providecommand{\href}[2]{#2}

\end{document}